\documentclass[12pt]{article} 
\usepackage[utf8x]{inputenc}
\usepackage{tikz} 
\usepackage[all]{xy}
\usepackage{authblk}
\usepackage[babel=true]{csquotes}
\usepackage{amsmath,amsthm, amssymb,amsfonts}
\usepackage[hmargin=1in,vmargin=1in]{geometry}
\numberwithin{equation}{subsection} 
\xyoption{arc}
\usepackage{eucal}
\usepackage{indentfirst}
\usepackage{mathrsfs}
\usepackage{verbatim}
\usepackage{makeidx}
\usepackage{graphicx}
\usepackage{yfonts} 
\usepackage{hyperref}
\usepackage{enumitem} 
\usepackage{extpfeil} 
\usepackage{dsfont}
\usepackage[T1]{fontenc}
\newtheorem{thm}{Theorem}[section]
\newtheorem{prop}[thm]{Proposition}

\newtheoremstyle{bidule}
{3pt}
{3pt}
{}
{}
{\scshape}
{.}
{.5em}
{}
\newtheorem{df}[thm]{Definition}

\theoremstyle{definition}

\newtheorem*{note}{Note}

\newtheorem*{warn}{Warning}

\newcommand{\C}{\mathcal{C}}
\newcommand{\Ub}{\mathcal{U}}

\newcommand{\F}{\mathcal{F}}

\newcommand{\Ar}{\text{Arr}}

\newcommand{\D}{\mathcal{D}}

\newcommand{\A}{\mathcal{A}}
\newcommand{\Aa}{\mathcal{A}}

\newcommand{\M}{\mathscr{M}}

\newcommand{\W}{\mathscr{W}}

\newcommand{\Un}{\mathbb{I}} 

\newcommand{\G}{\mathcal{G}}

\newcommand{\Cx}{\mathbb{C}}

\renewcommand{\to}{\longrightarrow}
\newcommand{\ol}{\overline}

\renewcommand{\k}{\textbf{k}} 

\newcommand{\tx}{\text}

\renewcommand{\to}{\longrightarrow}
\DeclareMathOperator\Id{Id}
\DeclareMathOperator\Hom{Hom}

\DeclareMathOperator\Set{\textbf{Set}} 
\DeclareMathOperator\Cat{\mathbf{Cat}}
\DeclareMathOperator\sCatb{\mathbf{sCat}_B}

\DeclareMathOperator\Ho{\mathbf{ho}} 

\DeclareMathOperator\Vect{\textbf{Vect}} 

\DeclareMathOperator\sset{\mathbf{sSet}} 
\DeclareMathOperator\ssetj{\mathbf{sSet}_{J}} 
\DeclareMathOperator\ssetq{\mathbf{sSet}_{Q}} 
\DeclareMathOperator\dCat{\mathbf{2-Cat}} 
 
\DeclareMathOperator\comsec{\M_\Ub[\ag]_{\tx{$e$}}^{\mathbf{c}}}

\DeclareMathOperator\comsepc{\M_\Ub[\ag]_{\tx{$e+$}}^{\mathbf{c}}}  
\DeclareMathOperator\com{\ag} 
 
\DeclareMathOperator\ag{\textgoth{A}}
\DeclareMathOperator\mua{\M_\Ub[\ag]}
\DeclareMathOperator\wua{\W_\Ub[\ag]}

\title{Understanding higher structures through Quillen-Segal objects} 
\author{Hugo V. Bacard \thanks{\textit{E-mail address}: \href{mailto:hbacard@uwo.ca}{hbacard@uwo.ca}
}}
 \affil{Western University}
\date{\today}
\begin{document}
\maketitle
\begin{abstract}
If $\M$ is a model category and $\Ub: \ag \to \M$ is a functor, we defined a Quillen-Segal $\Ub$-object as a weak equivalence $\F: s(\F) \xrightarrow{\sim} t(\F)$ such that $t(\F)=\Ub(b)$ for some $b\in \ag$.  
If  $\Ub$ is the nerve functor $\Ub: \Cat \to \ssetj$, with the Joyal model structure on $\sset$, then studying the comma category $(\ssetj \downarrow \Ub)$ leads naturally to concepts, such as Lurie's $\infty$-operad. It also gives simple examples of presentable, stable $\infty$-category, and higher topos. If we consider the \emph{coherent nerve} $\Ub: \sCatb \to \ssetj$,  then the theory of QS-objects directly connects with the program of Riehl and Verity.  If we apply our main result when $\Ub$ is the identity $\Id: \ssetq \to \ssetq$, with the Quillen model structure, the homotopy theory  of QS-objects is equivalent to that of Kan complexes and we believe that this is an \emph{avatar} of Voevodsky's \emph{Univalence axiom}. This equivalence holds for any combinatorial and left proper $\M$. This result agrees with our intuition, since by essence the `\emph{Quillen-Segal type}' is the \emph{Equivalence type}.
\end{abstract}
\setcounter{tocdepth}{1}

\vspace*{3cm}

\section{Background and Motivations}

This short discussion is motivated by our desire to have an understanding of the theory of quasicategories developed by Joyal and Lurie. And we also hope that our approach can be brought to a more general context of $(\infty,n)$-categories.\\

It has been proved by  Bergner \cite{Bergner_inf}, Joyal \cite{Joyal_quasi}, Lurie \cite{Lurie_HTT}, Rezk \cite{Rezk_mh},  that simplicial categories, quasicategories, complete Segal spaces and Segal categories, are all models for  $(\infty,1)$-categories. \\

There is an ongoing program developed by Riehl and Verity (see \cite{Riehl_Verity_adj, Riehl_Verity_2cat, Riehl_Verity_compl}), that aims to understand quasicategories through simplicial categories. And it turns out that this program coincide with the study of $\Ub$-QS-objects for the \emph{coherent nerve} of Cordier and Porter \cite{PC-hcc}:
$$\Ub: \sCatb  \to  \ssetj.$$

Here $\sCatb$ means that we consider the Bergner model structure on simplicial categories (see \cite{Bergner_scat}). 
Our discussion is based on  the previous paper  \cite{Bacard_QS}, wherein we expose the general idea of Quillen-Segal objects.\\ 

When we consider the comma category $(\ssetj\downarrow \Ub)$ we keep a control of what happens between simplicial categories and quasicategories. In other words we \emph{`temper' quasicategories} by linking them with rigid structures that are known to be equivalent. An object of this category is a morphism $\F: s(\F) \to t(\F)$, where $t(\F)= \Ub(\C)$ is a simplicial category. A morphism $\sigma: \F \to \G$ is given by a morphism of quasicategories $\sigma_0: s(\F) \to s(\G)$ and a morphism of simplicial categories $\sigma_1: \C \to \D$ such that we have a commutative square:
\[
\xy
(0,18)*+{s(\F)}="W";
(0,0)*+{t(\F)}="X";
(30,0)*+{t(\G)}="Y";
(30,18)*+{s(\G)}="E";
{\ar@{->}^-{\Ub(\sigma_1)}"X";"Y"};
{\ar@{->}^-{\F}"W";"X"};
{\ar@{->}^-{\sigma_0}"W";"E"};
{\ar@{->}^-{\G}"E";"Y"};
\endxy
\] 

We are not really interested in this big category, but rather the full subcategory $(\W_J \downarrow \Ub)$ of QS-objects $\F: s(\F) \xrightarrow{\sim} t(\F)$. We can define the weak equivalences as the morphism $\sigma: \F \to \G$ such that $\sigma_0$ is a weak equivalence in the Joyal model structure. The $3$-for-$2$ property will force $\sigma_1$ to be also a weak equivalence.

Now since $\Ub: \sCatb \to \ssetj$ is part of a Quillen equivalence it's not hard to show that the commutative triangle below descends to a triangle of equivalences between the homotopy categories. In fact we have a better statement as we shall see.
\[
\xy
(0,18)*+{(\W\downarrow \Ub)}="W";
(0,0)*+{\sCatb}="X";
(30,0)*+{\ssetj}="Y";
{\ar@{->}^-{\Ub}"X";"Y"};
{\ar@{->}^-{t}"W";"X"};
{\ar@{->}^-{s}"W";"Y"};
\endxy
\] 

With this equivalence of homotopy categories, we can pretend that the theory of quasicategories and simplicial categories can complete each other in the following sense. 
\begin{enumerate}
\item We can take advantage of the material developed by Joyal and Lurie for quasi-categories and transpose everything in term of simplicial categories, e.g $\infty$-limit and $\infty$-colimit, higher topos, etc.
\item Another relevant aspect is the result of Riehl-Verity \cite{Riehl_Verity_adj} about adjoint functor between quasicategories. Indeed, since we know that QS-objects also model $(\infty,1)$-categories, we have two two implication for a map  $\sigma: \F \to \G$  between QS-objects:
\begin{itemize}[label=$-$]
\item If the source $\sigma_0$ is an adjoint of quasicategories in the sense of Joyal-Lurie, then the target $\sigma_1$ should be a \emph{coherent adjoint} between simplicial categories.
\item Conversely if $\sigma_1$ is part of a homotopy coherent adjunction, then $\sigma_0$ should be an adjoint between quasicategories. 
\end{itemize} 
The result of Riehl-Verity says this indeed the case.
\end{enumerate} 

\paragraph*{Environment:} We fix a combinatorial and left proper model category $\M$ . We also fix a functor $\Ub: \ag \to \M$ that is a right adjoint between locally presentable categories. In \cite{Bacard_QS}, we also required $\Ub$ to be faithful, but this is not necessary. We need the faithfulness when we want to have an embedding $\ag \hookrightarrow (\M \downarrow \Ub)$. 

Let's recall the  definition of our objects of study.
\begin{df}
Let $\M$ and $\Ub$ be as above. 
\begin{enumerate}
\item Define the category $\mua$ of all $\Ub$-preobjects to be the comma category $(\M\downarrow \Ub)$. 
\item Define the category of Quillen-Segal $\Ub$-objects, henceforth $\Ub$-QS-objects,  to be the full subcategory $(\W \downarrow \Ub) \subseteq (\M\downarrow \Ub)$, where $\W \subseteq \M$ is the subcategory of weak equivalences.  We will denote this category by $\wua$.
\end{enumerate}
\end{df}

The main theorem in \cite{Bacard_QS} is the content of \cite[Theorem 8.2]{Bacard_QS}. We summarize the relevant information hereafter. 

\begin{thm}\label{main-thm-zip}
Let $\M$ be combinatorial and left proper model category. Let $\Ub: \com \to \M$ be a right adjoint between locally presentable categories.  Then the following hold. 
\begin{enumerate}
\item There are two combinatorial model structures on $\mua$, denoted by $\comsec$ and $\comsepc$. Both are left proper and the identity gives a Quillen equivalence:
$$ \comsec \to \comsepc.$$
\item In the  model category $\comsepc$, any fibrant object $\F: s(\F) \to \underbrace{t(\F)}_{=\Ub(b)}$ has the following properties.
\begin{itemize}[label=$-$]
\item $s(\F)$ is fibrant in $\M$;
\item The map $\F$ itself is a trivial fibration in $\M$. In particular, if $t(\F)$ is cofibrant in $\M$, then there is a generic weak \emph{`Voevodsky section'} $\pi: t(\F) \to s(\F)$ that is also a weak equivalence.
\end{itemize}
\item Assume that $\Ub: \com \to \M$ is a right Quillen functor that preserves and reflects the weak equivalences. Then the adjunction
$$ |-|^{\mathbf{c}}: \comsec \leftrightarrows \com: \iota,$$
is a Quillen equivalence.
\item In particular we have a diagram of equivalences between the homotopy categories.
$$\Ho[\comsepc] \xleftarrow{\simeq} \Ho[\comsec] \xrightarrow{\simeq} \Ho[\com].$$
\end{enumerate}
\end{thm}

\begin{note}
As outline before, in \cite[Theorem 8.2]{Bacard_QS}, we stated this theorem under the hypothesis that $\Ub$ is faithful but this is not necessary and is too restrictive. This hypothesis was suggested by the case of algebras and categories.  We can also show that if $\Ub$ is just part of a Quillen equivalence then we still have a Quillen equivalence between $ |-|^{\mathbf{c}}: \comsec \leftrightarrows \com: \iota$, and by closure, the other functor $\M \to \comsec$ is also part of a Quillen equivalence. We should warn the reader that if $\Ub$ is not faithful then the left adjoint $|-|$ is not the target functor anymore. 
\end{note}
We remind the reader that both model categories $\ssetq$ and $\ssetj$ are combinatorial and left proper.

\section{Applications and Interpretations}

\subsection{Understanding the Univalence axiom} 
The discussion that follows is an attempt to understand the univalence axiom, by a non-expert.\\
 
Let us consider the identity functor $\Ub=\Id: \ssetq \to \ssetq$ with the Quillen model structure. Then the category $\mua$ is just the category of morphisms in $\sset$: $$\Ar(\sset)=\sset^{\Un}= \Hom([0 \to 1], \sset)= Path(\sset).$$

The left adjoint $|-|:\sset^\Un \to \sset$ is the target-functor whose right adjoint is the embedding $\iota: \sset \to \sset^\Un$ that takes $X$ to the identity (type) morphism $\Id_X$. As far as we understand, the Quillen-Segal type in this context is precisely the equivalence type. \\

By the second assertion of Theorem \ref{main-thm-zip}, there exists a model structure on $\ssetq^\Un= Path(\ssetq)$ such that:
\begin{enumerate}
\item Every fibrant $\F$ object is a trivial fibration $\F: s(\F) \xtwoheadrightarrow{\sim} t(\F)$, with $s(\F)$ fibrant in $\ssetq$, which means that $s(\F)$ is Kan.


\item The underlying map $\F$ defines a tautological map $\F \to \Id_{t(\F)}$ in $\sset^\Un=Path(\sset)$. This map corresponds to the tautological commutative square:
\[
\xy
(0,18)*+{s(\F)}="W";
(0,0)*+{t(\F)}="X";
(30,0)*+{t(\F)}="Y";
(30,18)*+{t(\F)}="E";
{\ar@{->}^-{\Id_{t(\F)}}"X";"Y"};
{\ar@{->}^-{\F}"W";"X"};
{\ar@{->}^-{\F}"W";"E"};
{\ar@{->}^-{\Id_{t(\F)}}"E";"Y"};
\endxy
\]
\item Since every object in $\ssetq$ is cofibrant, then there is a \emph{Voevodsky section} $$\pi: t(\F) \to s(\F)$$ which is necessarily a weak equivalence. And by definition of a section we have an equality (proof):
 $$\pi \circ \F = \Id_{t(\F)}.$$
\item The section $\pi$ also determines a tautological map 
\begin{equation}\label{voev_sec}
\ol{\pi}: \Id_{t(\F)} \xrightarrow{\sim} \F,
\end{equation}
that corresponds to the commutative square:
  
\[
\xy
(0,18)*+{t(\F)}="W";
(0,0)*+{t(\F)}="X";
(30,0)*+{t(\F)}="Y";
(30,18)*+{s(\F)}="E";
{\ar@{->}^-{\Id_{t(\F)}}"X";"Y"};
{\ar@{->}^-{\Id_{t(\F)}}"W";"X"};
{\ar@{->}^-{\pi}"W";"E"};
{\ar@{->}^-{\F}"E";"Y"};
\endxy
\]
\item If we look at the respective type of the target and source of the  map $\ol{\pi}$, we would like to write something like:
$$(\quad = \quad ) \xrightarrow{\sim} ( \quad \simeq \quad).$$

Although we're not sure, it seems that this map somehow outlines the universality of the identity type which is the idea behind the Univalence axiom. If we abstract the Voevodsky section, then we are tempted to rewrite  \eqref{voev_sec} as:
$$V: Eq(\Id,QS).$$
\end{enumerate}

\subsection{Connection with the Riehl-Verity program}
As mentioned before there are two cases that are very interesting and that lead directly to the work of Riehl and Verity. 
\begin{enumerate}
\item If we consider the embedding $\Ub_1:\dCat \hookrightarrow \sCatb$ (see \cite{Riehl_Verity_2cat}). This leads somehow to  the  \emph{$2$-category theory of quasicategories } as named by the authors. 
\item And if we consider the coherent nerve $\Ub_2: \sCatb \to \ssetj$, this leads to \cite{Riehl_Verity_adj}.
\end{enumerate}
If we don't worry about size issues, the category $\sCatb$ is a combinatorial and left proper model category and each of $\Ub_1$ and $\Ub_2$ are right adjoint between locally presentable categories. Then by Theorem \ref{main-thm-zip} we have the following. 

\begin{prop}
\begin{enumerate}
\item There is a model structure on $(\sCatb \downarrow \Ub_1)$ in which any fibrant object is a trivial fibration of $\F: s(\F) \to \Ub_1(\C)$ with $s(\F)$ fibrant in $\sCatb$. In particular $s(\F)$ is enriched over Kan complexes and therefore the coherent nerve $\Ub_2(s(\F))$ is a quasicategory. 
\item There is a model structure on $(\ssetj \downarrow \Ub_2)$ in which any fibrant object is a trivial fibration of $\F: s(\F) \to \Ub_1(\C)$ with $s(\F)$ fibrant in $\ssetj$, that is a quasicategory.
\item The following chain of functor connects the two theories:
$$\dCat \xhookrightarrow{ \Ub_1} \sCatb \xrightarrow{\Ub_2} \ssetj.$$
\item Since in the Joyal model structure every object is cofibrant, there exists a section $\pi: \Ub_i(\C) \to s(\F)$ which is automatically an equivalence of quasicategories. Thanks to this section we can lift concepts of $2$-categories/simplicial categories to quasicategories and vice-versa.
\end{enumerate}
\end{prop}

\begin{warn}
There is a potential conflict of terminology between ours and that of Riehl-Verity. If we let $\M$ be either $\ssetj$ or $\sCatb$, then with our notation,  the two comma categories will be denoted by $\M_{\Ub_1}[\dCat]$ and $\M_{\Ub_2}[\sCatb]$. In \cite{Bacard_QS}, we've interpreted the notation $\mua$ as the category of objects of $\M$ with \emph{coefficients or coordinates} in $\ag$ (with respect to $\Ub$).

In particular we would say for example that: ``$\M_{\Ub}[\dCat]$ is the quasicategory theory of $2$-categories !''. This terminology is justified by two reasons:
\begin{itemize}[label=$-$]
\item The first reason is because the fibrant objects there are the quasicategories that are equivalent to $2$-categories;
\item The second reason is the  language of \emph{points} introduced by Grothendieck, as we consider the representable $\Hom(-,\M)$. Moreover it's highly likely that the forgetful functor we've considered are geometric morphisms of higher topoi.
\end{itemize}
\end{warn}

\subsection{Producing examples of `higher concepts'}
In this part, we simply list some examples that follow from  Theorem \ref{main-thm-zip} applied to the usual nerve functor 
$$\Ub: \Cat \to \ssetj.$$ 

In this case the homotopy theory $\mua$ is the \emph{Quasicategory theory of $1$-categories}. Indeed, just like previously,  any fibrant object is a trivial fibration $\F: s(\F) \to \Ub(\C)$, such that $s(\F)$ is fibrant i.e., a quasicategory and  where $\C$ is a usual $1$-category.\\

As usual, since we have a trivial fibration whose source is fibrant and the target is cofibrant, there is a section $\pi: \C \to s(\F)$. Given such section, it's an interesting exercise to determine `by hands', the structure of the quasicategory $s(\F)$ attached to $\C$, if:
\begin{enumerate}
\item $\C$ is usual abelian category;
\item $\C$ is the derived category of an abelian category;
\item $\C$ is (pre)triangulated
\item $\C$ is a Grothendieck topos, 
\item $\C$ is complete/cocomplete;
\item $\C$ is the stable homotopy category,
\item $\C$ is a symmetric monoidal category, Tannakian;
\item $\C$ is a (stable) model category, combinatorial model category and so on.
\end{enumerate}

We know that we should land with examples of the notions introduced by Joyal \cite{Joyal_quasi} and Lurie \cite{Lurie_HALG,Lurie_HTT}. \\

Even if we don't consider the fibrant objects, but just objects $\F: s(\F) \to \Ub(\C)$ satisfying some conditions, we end up with concepts such as Lurie's $\infty$-operads.\\

The category of fibrant objects is the category of quasicategories that are the most close to $1$-categories. For example if we have any functor of $1$-category such as a localization functor going to the derived category:
$$L: \Aa \to D(\Aa),$$  then a fibrant replacement of $D(\Aa)$ is a trivial fibration $\F: s(\F) \to \Ub(D(\Aa))$. And because $\Aa$ is cofibrant in the Joyal model structure, we can find a lift  $\Aa \to s(\F)$. This means that $s(\F)$ is the best quasicategory that approaches the derived category $D(\Aa)$.  In fact, homing a trivial fibration between quasicategories yields an equivalence of quasicategories.

\subsection{General picture: Minimal Homotopy Enhancement}
There are many categories such as $\Set$ or  $\Vect_\k$, that  don't have an interesting homotopy theory.  The general idea of Quillen-Segal objects is precisely to enhance these category by embedding them into a better category where there is a good homotopy theory: this is what we call \emph{homotopy enhancement}.

Classically we would just embed $\Ub:\Set \hookrightarrow \sset$ and similarly  we would embed $\Ub: \Vect_\k \hookrightarrow C(\k)$ or to simplicial vector spaces. But these embedding are too big. Moreover, there are many reasons why we would like to embed them directly to quasicategories or to Kan complexes, rather than passing for example through Dold-Kan, etc.\\

With the Quillen-Segal formalism we get a \emph{minimal embedding} each time. Here is how we would interpret these embeddings. 
\begin{enumerate}
\item The embedding $\Vect_\k \hookrightarrow (\ssetj \downarrow \Ub)$ gives the quasicategory theory of vector spaces. Any fibrant object there is a quasicategory $s(\F)$ that looks like a vector space ! 
\item We can alternatively embed $\Vect_\k \hookrightarrow (\ssetq \downarrow \Ub)$, and this should be the right direction toward type theory (Kan complexes).
\item We also have two possibilities $\Set \hookrightarrow (\ssetq \downarrow \Ub)$ or $\Set \hookrightarrow (\ssetq \downarrow \Ub)$. 
\item If $\M= \ssetj$, and $\ag=\Set$, we get the quasicategory theory of sets. Here we regard the inclusion $\Set \to \ssetj$ as a right adjoint, and indeed as geometric point of the higher topos $\ssetj$. In particular it's a monoidal functor with the cartesian product. 

We can put the relative pushout product on $\mua$. This way we get a monoidal category of quasicategory of sets. In particular we can determine the controlled quasicategories $s(\F) \to X$ that look like:
\begin{itemize}[label=$-$]
\item The integers,
\item The rationals,
\item The real numbers;
\item The complex numbers;
\item a ring, field, groups, Hodge structure (...and maybe a prime number in a quasicategory)
\end{itemize}

It's reasonable to believe that this should join Lurie's program on \emph{Higher Algebra} \cite{Lurie_HALG}. Morally we're enhancing the Grothendieck topos $(\Set, \times, 1)$ by a \emph{minimal homotopical model}. And the advantage of this is that we keep track of everything since there is always an weak section (Voevodsky section) $\pi: X \to s(\F)$.
\item Doing the same thing with the Quillen model structure $\ssetq$ should normally fits in type theory.

\item We can restrict to the right adjoint $\Vect_\k \to \Set \to \ssetj$, and do the relative lax pushout product. The category $\mua$, that is the quasicategory of vectors spaces will inherit a monoidal product where the unit is the map $\ast \to \k$ that selects the unit element of the field $\k$. We still have a functor that is monoidal $\mua \to \Vect_\k$.\\

Doing this is not only general nonsense. We can use for example the category $\ssetj[\Vect_\Cx]$ as coefficient category for general TQFT and so on.
\item The category $(\Set,\times,1)$ is known as the Grothendieck topos of sheaves over the point. Given a general site $\C$, we can form a homotopical minimal enhancement of its category of sheaves by taking $\M$ to be the model category of simplicial (pre)sheaves \emph{à la} Jardine-Joyal. The category $\mua$ in this case is still a combinatorial and left proper and we can perform again Bousfield localization to \emph{converge} to Morel-Voevodsky work \cite{Mo_Vo_A1}.
\item Finally it's important to observe that  if $\Ub:\A \hookrightarrow Mod(\A^{op}$ is the Yoneda embedding of a dg-category, then  a QS-object is exactly what To\"{e}n call quasi-representable (see \cite{Toen_Morita}).
\end{enumerate} 

\begin{center}
Keep enhancing !
\end{center}

\bibliographystyle{plain}
\bibliography{Bibliography_LP_COSEG}
\end{document}